\documentclass{amsart}

\title   [Hopf algebras in positive characteristic]
         {Semisimplicity criteria for irreducible Hopf algebras in positive characteristic}
\author  {Akira Masuoka}
\address {Institute of Mathematics, University of Tsukuba, Ibaraki 305-8571, Japan}
\email   {akira@math.tsukuba.ac.jp}
\date    {}

\subjclass[2000]{16W30}
\keywords{Hopf algebra in positive characteristic, restricted Lie algebra}

\usepackage{amsmath}
\usepackage{amssymb}
\usepackage{amsthm}

\theoremstyle{plain}
\newtheorem{theorem}                 {Theorem}[section]
\newtheorem{lemma}          [theorem]{Lemma}
\newtheorem{proposition}    [theorem]{Proposition}

\theoremstyle{definition}
\newtheorem{definition}     [theorem]{Definition}

\numberwithin{equation}{section}

\theoremstyle{remark}
\newtheorem{remark}         [theorem]{Remark}

\newcommand{\id}            {\mathrm{id}}

\newcommand{\hit}           {\rightharpoonup}
\newcommand{\ch}            {\mathop{\rm ch}\nolimits}
\newcommand{\gr}            {\mathop{\rm gr}\nolimits}
\newcommand{\Ker}           {\mathop{\rm Ker}\nolimits}

\begin{document}

\begin{abstract}
  We prove that a finite-dimensional irreducible Hopf algebra $H$ in positive characteristic is semisimple,
  if and only if it is commutative and semisimple,
  if and only if the restricted Lie algebra $P(H)$ of the primitives is a torus.
  This generalizes Hochschild's theorem on restricted Lie algebras,
  and also generalizes Demazure and Gabriel's and Sweedler's results on group schemes,
  in the special but essential situation with finiteness assumption added.
\end{abstract}

\maketitle

\setcounter{section}{-1}

\section{Introduction}

During the last two decades or so,
much progress has been brought to the study of semisimple Hopf algebras; see \cite{bib:M}, for example.
Let $H$ be a finite-dimensional Hopf algebra over a field $k$ with antipode $S$.
Let us consider the following conditions.
\begin{enumerate}
\renewcommand{\labelenumi}{(\hbox to 1.25ex {\hfill{}\alph{enumi}\hfill})}
\item $H$ is semisimple as an algebra.
\item $H$ is cosemisimple as a coalgebra.
\item The characteristic $\ch k$ does not divide the dimension $\dim H$.
\item $S$ is an involution, i.e., $S \circ S = \id_H$.
\end{enumerate}
Summarizing results by Larson \cite{bib:L}, Larson and Radford \cite{bib:LR},
and Etingof and Gelaki \cite{bib:EG}, we have that (a) \& (b) $\iff$ (c) \& (d),
and that if $\ch k = 0$ or $> \sqrt{d}^{\varphi(d)}$, where $d = \dim H$, then (a) $\iff$ (b) $\iff$ (d).
Moreover, Etingof and Gelaki \cite{bib:EG} showed that by `lifting',
most of known results on semisimple (then necessarily cosemisimple) Hopf algebras in characteristic zero
can extend to semisimple and cosemisimple Hopf algebras in positive characteristic.
However, little seems known on those semisimple Hopf algebras in positive characteristic which are not cosemisimple.
In this paper we will determine when a finite-dimensional irreducible Hopf algebras in positive characteristic is semisimple.
Recall from \cite[Sect. 8.0]{bib:Sw1}
that a Hopf algebra $H$ over a field $k$ is said to be {\em irreducible} (as a coalgebra),
if the coradical of $H$ equals $k$,
or equivalently if the trivial $H$-comodule $k$ is a unique (up to isomorphism) simple $H$-comodule;
this property is very opposite to cosemisimplicity.
In this case the dimension $\dim H$, if finite, equals a power of $p$,
provided $\ch k = p > 0$; see Proposition \ref{prop:1.1} (1).

Throughout we work over a field $k$.
Let $\bar{k}$ denote the algebraic closure of $k$.
Given a finite group $G$, we let $k^G$ denote the dual (commutative semisimple) Hopf algebra of the group algebra $kG$.
Our aim is to prove the following.

\begin{theorem} \label{thm:0.1}
  Suppose $\ch k = p > 0$.
  Let $H$ be an irreducible Hopf algebra of finite dimension.
  Then the following are equivalent.

  \begin{itemize}
  \item [(i)]   $H$ is semisimple.
  \item [(ii)]  $H$ is commutative and semisimple.
  \item [(iii)] $H \otimes \bar{k} \cong \bar{k}^G$, where $G$ is a finite $p$-group.
  \item [(iv)]  The restricted Lie algebra $P(H)$ of all primitives in $H$ is a torus; see Definition \ref{def:1.2}.
  \end{itemize}
\end{theorem}

This generalizes Hochschild's theorem on restricted Lie algebras,
and also generalizes Demazure and Gabriel's and Sweedler's results on group schemes,
in the special but essential situation with finiteness assumption added.
Keep the assumption $\ch k = p > 0$.
For a restricted Lie algebra $L$ of finite dimension,
the restricted enveloping algebra $u(L)$ is a cocommutative irreducible Hopf algebra of finite dimension.
Hochschild \cite{bib:H} proved the theorem above, when $H = u(L)$.
Supposing $k = \bar{k}$, Demazure and Gabriel \cite[IV, \S 3, 3.7]{bib:DG} proved that
a connected algebraic affine group scheme $\mathfrak{G}$ is diagonalizable,
if it does not have any closed subgroup scheme isomorphic to $\boldsymbol{\alpha}_p$.
Sweedler \cite{bib:Sw2} proved that a connected reductive affine group scheme $\mathfrak{G}$
is necessary abelian; see also \cite[IV, \S 3, 3.6]{bib:DG} referred to as Th\'eor\`eme de Nagata.
In the situation above, suppose that $\mathfrak{G}$ is finite,
or in other words that the coordinate algebra $\mathcal{O}(\mathfrak{G})$ is finite-dimensional.
Then the dual Hopf algebra $\mathcal{O}(\mathfrak{G})^*$ is cocommutative,
and irreducible by the connectedness assumption,
whence our theorem can apply to $\mathcal{O}(\mathfrak{G})^*$.
More precisely, the last two results with $\mathfrak{G}$ supposed to be finite are equivalent to our results,
$\mbox{(iv)} \Rightarrow \mbox{(iii)}$, $\mbox{(i)} \Rightarrow \mbox{(ii)}$, respectively,
with $H$ supposed to be cocommutative.
It is easy to deduce those two results on $\mathfrak{G}$ in general, as Sweedler \cite{bib:Sw2} actually did,
from the results in the special case of $\mathfrak{G}$ being finite; see Remark \ref{rem:3.2}.
The author does not yet know how our theorem can extend to infinite-dimensional Hopf algebras in the dual situation.

Our proof of the theorem is quite simple,
based on an elementary observation on what we call {\em relative primitives};
see the proof of our key lemma, Lemma \ref{lem:2.2}.
It does not depend on any of the three known results cited above, and differs from their proofs.

\section{Preliminaries}

Let $H$ be a Hopf algebras.
The coalgebra structures and the antipode are denoted by
\[
\Delta: H \to H \otimes H,
\quad \Delta(h) = \sum h_1 \otimes h_2,
\quad \varepsilon: H \to k,
\quad S: H \to H,
\]
respectively.
We let $H^+ = \Ker \varepsilon$ denote the augmentation ideal.

Suppose $\dim H < \infty$.
Then, $H$ is semisimple if and only if $H \otimes \bar{k}$ is semisimple.
Recall also that $H$ is commutative and semisimple if and only if $H \otimes \bar{k} \cong \bar{k}^G$,
where $G$ is a finite group.

The subspace
\[ P(H) = \{ x \in H \mid \Delta(x) = x \otimes 1 + 1 \otimes x \} \]
in $H$ consisting of all primitives forms a Lie algebra
with respect to the commutator $[x, y] = xy - yx$. We have
\begin{equation}
  \label{eq:1.1} P(H) \otimes \bar{k} = P(H \otimes \bar{k})
\end{equation}

Suppose that $H$ is irreducible.
Then the coradical filtration $k = H_0 \subset H_1 \subset \cdots$ makes $H$ into a filtered Hopf algebra
$H = \bigcup_{n \ge 0} H_n$, so that the associated graded Hopf algebra
$\gr H = \bigoplus_{n \ge 0} H_n / H_{n - 1}$ ($H_{-1} = 0$) is strictly graded in the sense that
$\gr H(0) = k$, $\gr H(1) = P(\gr H)$; see \cite[Sect. 11.2]{bib:Sw1}.
It is known that $\gr H$ is commutative \cite[Theorem 11.2.5]{bib:Sw1}.

In what follows until the end of this section,
we suppose $\ch k = p > 0$.
Given a restricted Lie algebra $L$,
the restricted enveloping algebra $u(L)$ naturally forms a cocommutative irreducible Hopf algebra.
It is characterized as such a Hopf algebra $H$ that is generated by $P(H)$, and satisfies $P(H) = L$.
If $\dim L = n$, then $\dim u(L) = p^n$.
In every Hopf algebra $H$, $P(H)$ forms a restricted Lie algebra with respect to $x^{[p]} = x^p$,
the $p^{\mbox{th}}$ power in $H$, and generates a Hopf subalgebra isomorphic to $u(P(H))$.

\begin{proposition} \label{prop:1.1}
  Let $H$ be an irreducible Hopf algebra of finite dimension.

  (1) Set $L = P((\gr H)^*)$ in the dual graded Hopf algebra $(\gr H)^*$ of $\gr H$.
  Then, $L$ is a positively graded restricted Lie algebra $\bigoplus_{n > 0} L(n)$ generated by $L(1)$,
  and $u(L) = (\gr H)^*$. In particular, $\dim H = p^n$, if $n = \dim L$.

  (2) $H$ is commutative and semisimple if and only if $H \otimes \bar{k} \cong \bar{k}^G$,
  where $G$ is a finite $p$-group.

  (3) $H$ is commutative, semisimple and generated by $P(H)$,
  if and only if $u(P(H)) = H$ and $P(H)$ is a torus(see Definition \ref{def:1.2} below),
  if and only if $H \otimes \bar{k} \cong \bar{k}^{G}$, where $G$ is a $p$-torsion finite abelian group.
\end{proposition}

\begin{proof}
  (1) This follows since $(\gr H)^*$ is generated by the first component,
  and is cocommutative by \cite[Theorem 11.2.5]{bib:Sw1}.

  (2) This follows, since a finite group algebra $kG$ is local
  if and only if $G$ is a $p$-group.

  (3) Suppose that $G$ is a finite $p$-group.
  Notice then that $\bar{k}^G$ is generated by primitives,
  if and only if $\bar{k}G$ is commutative and $x^p = 0$ for every element $x$ in the augmentation ideal of $\bar{k}G$,
  if and only if $G$ is abelian and $p$-torsion.
  The desired result now follows from (2) above and the next proposition, (a) $\iff$ (e).
\end{proof}

\begin{definition}[{\cite[p.86]{bib:SF}}] \label{def:1.2}
  Let $L$ be a restricted Lie algebra of finite dimension.
  $L$ is called a {\em torus}, if
  \begin{itemize}
  \item [(1)] $L$ is abelian, and
  \item [(2)] every element of $L$ is semisimple in $u(L)$, i.e., generates a semisimple subalgebra.
  \end{itemize}
\end{definition}

\begin{proposition}[Hochschild \cite{bib:H}] \label{prop:1.3}
  For a restricted Lie algebra $L$ of finite dimension, the following are equivalent.
  \begin{enumerate}
    \renewcommand{\labelenumi}{(\hbox to 1.25ex {\hfill{}\alph{enumi}\hfill})}
  \item $L$ is a torus.
  \item $L$ is abelian, and is spanned by $L^{[p]}$.
  \item $L \otimes \bar{k}$ does not contain any non-zero element $z$ with $z^{[p]} = 0$.
  \item $u(L)$ is semisimple.
  \item $u(L)$ is commutative and semisimple.
  \end{enumerate}
\end{proposition}

\begin{remark} \label{rem:1.4}
  To prove Theorem \ref{thm:0.1}, we will not use this last proposition.
  On the contrary, the equivalences (c) $\iff$ (d) $\iff$ (e) follow from (the proof of) the theorem;
  see Remark \ref{rem:3.1}.
  In particular, (c) (or (d)) implies that $L$ is abelian.
  Let us record here a proof of the remaining.

  (e) $\Longrightarrow$ (a).
  In $u(L)$, the subalgebra $k[x]$ generated by an element $x \in L$ is a Hopf subalgebra,
  which is semisimple if $u(L)$ is; see \cite[3.2.3, p.31]{bib:Mo}.

  (a) $\Longrightarrow$ (e).
  In general, the minimal polynomial of $0 \ne x \in L$ in $u(L)$ is of the form
  \begin{equation} \label{eq:1.2}
    c_0 x + c_1 x^p + \cdots + c_{r-1} x^{p^{r-1}} + x^{p^r} \quad (r > 0, c_i \in k),
  \end{equation}
  where $x, x^p, \cdots, x^{p^r}$ are linearly dependent in $L$ with $r$ minimal.
  The element $x$ is semisimple if and only if $c_0 \ne 0$,
  in which case $x$ generates a separable subalgebra.
  Condition (a) implies that in the affine algebra $u(L)$, the largest separable subalgebra equals to $u(L)$.

  (b) $\iff$ (c).
  We may assume that $L$ is abelian.
  Let $k^{1/p}$ denote the subfield of $\bar{k}$ (including $k$)
  consisting of the $p^{\mbox{th}}$ roots of all elements in $k$.
  By the last assumption we have a map
  \[ \phi_L: L \otimes k^{1/p} \to L, \quad \phi_L(x \otimes c^{1/p}) = c x^{[p]}. \]
  Regard $L \otimes k^{1/p}$ as a $k$-vector space through $k \to k^{1/p}$, $c \mapsto c^{1/p}$.
  Then, $\phi_L$ is $k$-linear, and the base extension $\phi_L \otimes \bar{k}$ to $\bar{k}$ is identified with
  $\phi_{L \otimes \bar{k}}$.
  Therefore, (b) and (c) are both equivalent to the bijectivity of $\phi_L$.
\end{remark}

\section{The space $P(H, K)$ of relative primitives}

Before Lemma \ref{lem:2.2}, the characteristic $\ch k$ may be arbitrary.
Let $H$ be a Hopf algebra, and let $K \subset H$ be a Hopf subalgebra.
We define a subspace of $H$ by
\[ P(H, K) = \{ x \in H^+ \mid \Delta(x) - x \otimes 1 - 1 \otimes x \in K \otimes K \}, \]
whose elements may be called {\em relative primitives}. We see $K^+ \subset P(H, K)$. Let
\begin{equation} \label{eq:2.1}
  a \hit h = \sum a_1 h S(a_2) \quad (a, h \in H)
\end{equation}
denote the conjugation by $H$ on itself.
Then, $H$ is a left $H$-module under $\hit$.

\begin{proposition} \label{prop:2.1}
  (1) If $B \subset K$ is a cocommutative Hopf subalgebra,
  $P(H, K)$ is a $B$-submodule of $H$ under the conjugation.

  (2) Let $k\langle P(H, K) \rangle$ denote the subalgebra of $H$ generated by $P(H, K)$.
  Then this is a Hopf subalgebra including $K$.

  (3) Suppose that $H$ is irreducible.
  If $H \supsetneq K$, then $P(H, K) \supsetneq K^+$.
\end{proposition}

\begin{proof}
  (1) Since in general,
  \begin{equation} \label{eq:2.2}
    \Delta(a \hit h) = \sum a_1 h_1 S(a_3) \otimes (a_2 \hit h_2),
  \end{equation}
  it follows that if $a \in B$, then $\Delta(a \hit h) = \sum (a_1 \hit h_1) \otimes (a_2 \hit h_2)$.
  This implies the assertion.

  (2) This follows since one sees that the sum $k + P(H, K)$ forms a subcoalgebra of $H$ stable under the antipode.

  (3) Since $H$ is irreducible, the filtered Hopf subalgebra $\bigcup_{n \ge 0} F_n H$ ($\subset H$)
  given by the wedge products $F_n H = \bigwedge^{n + 1} K$ coincides with $H$; see \cite[Sect. 5.2]{bib:Mo}.
  Explicitly, we have $F_0 H = K$, and
  \[ F_1 H = \{ h \in H \mid \Delta(h) \in K \otimes H + H \otimes K \}. \]
  Let $\mathcal{H} = \bigoplus_{n \ge 0} F_n H / F_{n - 1} H$ ($F_{-1} H = 0$)
  denote the associated graded Hopf algebra.
  As is easily seen, $F_1 H$ is a subcoalgebra of $H$,
  and the natural $(F_1 H, F_1 H)$-bicomodule structure on $F_1 H$ induces a $(K, K)$-bicomodule structure,
  say $K \otimes \mathcal{H}(1) {\displaystyle \overset{\rho_L}{\longleftarrow}} \mathcal{H}(1)
  {\displaystyle \overset{\rho_R}{\longrightarrow}} \mathcal{H}(1) \otimes K$, on $\mathcal{H}(1)$
  (this holds for all $F_n H$, $\mathcal{H}(n)$).
  Let \[ \mathcal{H}(1)_0 = \{ u \in \mathcal{H}(1) \mid \rho_L(u) = 1 \otimes u, \ \rho_R(u) = u \otimes 1 \} \]
  denote the subspace in $\mathcal{H}(1)$ consisting of all coinvariants on both sides.
  We see that $k + P(H, K)$ is included in $F_1 H$, and coincides with the pullback of $\mathcal{H}(1)_0$
  along the natural projection $F_1 H \to \mathcal{H}(1)$.
  Hence we have a linear surjection $P(H, K) \to \mathcal{H}(1)_0$, whose kernel equals $K^+$,
  and which therefore induces an isomorphism
  \begin{equation} \label{eq:2.3}
    P(H, K) / K^+ \cong \mathcal{H}(1)_0 .
  \end{equation}
  The assumption $K \subsetneq H$ implies $K \subsetneq \mathcal{H}$,
  which in turn implies $\mathcal{H}(1) \ne 0$ since the iterated coproduct induces a monomorphism
  $\mathcal{H}(n) \to \mathcal{H}(1)^{\otimes n}$ for each $n > 0$.
  Since $K$ is irreducible, $\mathcal{H}(1)_0$ coincides with the socle of the $(K, K)$-bicomodule $\mathcal{H}(1)$,
  and hence is non-zero. By (\ref{eq:2.3}), this proves $P(H, K) \supsetneq K^+$, as desired.
\end{proof}

Keep the situation as above.
The following is a key result to prove Theorem \ref{thm:0.1}.

\begin{lemma} \label{lem:2.2}
  Suppose $\ch k = p > 0$, and that $k$ is perfect.
  Suppose $\dim H < \infty$, and $K \cong k^G$, where $G$ is a finite $p$-group.
  Suppose in addition that $H$ does not contain any non-zero primitive $z$ with $z^p = 0$.
  Then, $K$ is included in the center of $k \langle P(H, K) \rangle$.
\end{lemma}

\begin{proof}
  Being a finite $p$-group, $G$ has a central series,
  $G = G_0 \supset G_1 \supset \cdots \supset G_s = \{ 1 \}$.
  Thus for each $0 < i \le s$, $G_i$ is normal in $G$, and $G_{i-1} / G_i$ is central in $G/G_i$.
  Set $\bar{G}_i = G / G_i$, $B_i = k^{\bar{G}_i}$. Then we have inclusions
  \[ k = B_0 \subset B_1 \subset \cdots B_s = K \]
  of commutative Hopf algebras.
  Set \[ \bar{B}_i = k^{G_{i-1} / G_i} ( = B_i / B_iB_{i-1}^+). \]
  Then the Hopf algebra quotient $B_i \to \bar{B}_i$, $b \mapsto \bar{b}$ is cocentral in the sense
  \begin{equation} \label{eq:2.4}
    \sum b_1 \otimes \bar{b}_2 = \sum b_2 \otimes \bar{b}_1 \quad (b \in B_i).
  \end{equation}

  To prove the lemma we need to prove that $K$ trivially acts on $P(H, K)$ under the conjugation $\hit$ (see (\ref{eq:2.1})),
  or explicitly that $a \hit x = \varepsilon(a) x$ for all $a \in K$, $x \in P(H, K)$.
  The first step is to prove that $B_1$ trivially acts on $P(H, K)$.
  Since $K$ is commutative, $B_1$ (and moreover $K$) trivially acts on $K$.
  Since $B_1$ is cocommutative, it follows by Proposition \ref{prop:2.1} (1) that $P(H, K)$ is a $B_1$-submodule of $H$.
  Let $x \in P(H, K)$, and define
  \begin{equation} \label{eq:2.5}
    \gamma_x := \Delta(x) - x \otimes 1 - 1 \otimes x \quad (\in K \otimes K).
  \end{equation}
  Notice that
  \begin{equation} \label{eq:2.6}
    \gamma_{a \hit x} = \varepsilon(a) \gamma_x \quad (a \in B_1).
  \end{equation}
  Let $e_g$ ($g \in \bar{G}_1$) denote the primitive idempotents in $B_1 = k^{\bar{G}_1}$,
  so that $e_g(f) = \delta_{f, g}$ ($f \in \bar{G}_1$). We wish to prove
  \[ e_g \hit x = \delta_{1, g} x \quad (g \in \bar{G}_1). \]
  We may re-choose $x$ so that $B_1 \hit x$ is 1-dimensional,
  or in other words, so that $x \ne 0$, and
  \begin{equation} \label{eq:2.7}
    e_g \hit x = \delta_{t, g} x \quad (g \in \bar{G}_1)
  \end{equation}
  for some $t \in \bar{G}_1$.
  We aim to prove $t = 1$. Suppose $t \ne 1$, on the contrary.
  We then see that $x$ is a primitive, or $\gamma_x = 0$.
  In fact we see from (\ref{eq:2.6}), (\ref{eq:2.7}) that
  \[ \delta_{t, g} \gamma_x = \varepsilon(e_g) \gamma_x \ ( = \delta_{1, g} \gamma_x) \quad (g \in \bar{G}_1), \]
  which tells $\gamma_x = 0$ when $g = t$ ($\ne 1$).
  Notice now $x^{p^i} \in P(H)$ for each $i \ge 0$.
  Since $\dim H < \infty$, the subalgebra $k[x]$ of $H$ generated by $x$ has a defining relation of the form
  \begin{equation} \label{eq:2.8}
    c_0 x + c_1 x^p + \cdots + c_{r - 1} x^{p^{r-1}} + x^{p^r} = 0
    \quad (r > 0, c_i \in k).
  \end{equation}
  See (\ref{eq:1.2}). We have $c_i^{1/p}$ in $k$, since $k$ is perfect.
  We see $c_0 \ne 0$, since otherwise,
  \begin{equation} \label{eq:2.9}
    z := c_1^{1/p} x + \cdots + c_{r - 1}^{1/p} x^{p^{r-2}} + x^{p^{r-1}}
  \end{equation}
  would be a non-zero primitive with $z^p = 0$.
  Notice $e_g \hit x^{p^i} = \delta_{t^{p^i}, g} x^{p^i}$.
  By applying $e_g \hit$ to (\ref{eq:2.8}), we obtain
  \begin{equation} \label{eq:2.10}
    \delta_{t, g} c_0 x + \cdots + \delta_{t^{p^r}, g} x^{p^r} = 0.
  \end{equation}
  Choose $t^{p^r}$ as $g$.
  By comparing (\ref{eq:2.8}) with (\ref{eq:2.10}),
  it then follows that $t = g = t^{p^r}$, which implies t = 1, as desired.

  As the second step we wish to prove that $B_2$ trivially acts on $P(H, K)$.
  Let $a \in B_2$, $x \in P(H, K)$. Then, $a \hit x \in H$.
  By the result of the first step above,
  this last action factors through $\bar{B_2}$ so that $a \hit x = \bar{a} \hit x$.
  By (\ref{eq:2.2}) and (\ref{eq:2.4}), we see
  \begin{align*}
    \varepsilon(a) \gamma_x
    & = \Delta(a \hit x) - (a \hit x) \otimes 1 - \sum a_1 S(a_3) \otimes (\bar{a}_2 \hit x) \\
    & = \Delta(a \hit x) - (a \hit x) \otimes 1 - 1 \otimes (a \hit x)
  \end{align*}
  Therefore, $a \hit x \in P(H, K)$, and $P(H, K)$ is a $B_2$-submodule of $H$.
  In addition, (\ref{eq:2.6}) holds for $a \in B_2$, as well.
  The same argument as in the first step, but $\bar{G}_1$ replaced with $\bar{G}_2$,
  proves the desired result of this step.

  Putting forward the steps we see finally that $B_s ( = K)$ trivially acts on $P(H, K)$, as desired.
\end{proof}

\section{Proof of Theorem \ref{thm:0.1}}

Suppose that we are in the situation of the theorem.
The equivalence (ii) $\iff$ (iii) follows by Proposition \ref{prop:1.1} (1).
Obviously, (ii) $\Longrightarrow$ (i).
We see (ii) $\Longrightarrow$ (iv), as in the proof of (e) $\Longrightarrow$ (a) in Remark \ref{rem:1.4}.
To complete the proof we will see (i) $\Longrightarrow$ (ii), (iv) $\Longrightarrow$ (ii).

By base extension we may suppose $k = \bar{k}$.
Assume (i) or (iv).
Whichever is assumed, $H$ does not contain any non-zero primitive $z$ with $z^p = 0$,
since every primitive then needs to be semisimple, by \cite[3.2.3]{bib:Mo}; as for (iv),
see Remark \ref{rem:1.4}, (a) $\Longrightarrow$ (e).
Let $K \subset H$ be a commutative semisimple Hopf subalgebra of the largest dimension.
Then, $K \cong k^G$, where $G$ is a finite $p$-group,
and Lemma \ref{lem:2.2} can apply to $K \subset H$.
We wish to prove $K = H$.
Suppose $K \subsetneq H$, on the contrary.
By Proposition \ref{prop:1.3} (3), we have an element $x$ in $P(H, K) \setminus K$.
Then, $K$ and $x$ generate a sub-bialgebra, necessarily a Hopf subalgebra, say $J$, in $H$.
By Lemma \ref{lem:2.2}, $J$ is a (commutative) $K$-algebra generated by $x$, properly including $K$.
Keep in mind that by the commutativity,
the formula $(\alpha + \beta)^{p^i} = \alpha^{p^i} + \beta^{p^i}$ holds in $J$, and in $J \otimes J$.
We wish to prove that $J$ is semisimple, which will conclude $K = H$, as desired.

Let $Q = J / K^+ J$ denote the quotient Hopf algebra of $J$ divided by $K^+ J$.
Let $\pi: J \to Q$ denote the quotient map.
We regard $J$ as a right $Q$-comodule algebra with respect to $(\id \otimes \pi) \circ \Delta$.
Then the $Q$-coinvariants in $J$ are precisely $K$, that is,
\[ K = \{ h \in J \mid (\id \otimes \pi) \circ \Delta(h) = h \otimes 1 \}. \]
Set $y = \pi(x)$. Then, $y$ is a non-zero primitive, and generates $Q$.
A defining relation of $Q$ is given by an equation such as
\[ c_0 y + c_1 y^p + \cdots + c_{r-1} y^{p^{r-1}} + y^{p^r} = 0 \quad (r > 0, c_i \in k). \]
See (\ref{eq:1.2}). Define
\[ a := c_0 x + c_1 x^p + \cdots + c_{r-1} x^{p^{r-1}} + x^{p^r} \quad (\in H). \]
Then, $\pi(a) = 0$, and $a \in P(H, K)$ since now $x^{p^i} \in P(H, K)$.
Thus, $\Delta(a) \in a \otimes 1 + 1 \otimes a + K^+ \otimes K^+$.
It follows that $a$ is right $Q$-coinvariant, whence $a \in K$.
We may suppose $a = 0$, by replacing $x$ with $x - b$,
where $b \in K^+$ such that $c_0 b + c_1 b^p + \cdots + b^{p^r} = a$.
It then follows by \cite[7.2.2, p.106]{bib:Mo}
that $y \mapsto x$ defines a $Q$-colinear $K$-algebra isomorphism $K \otimes Q \cong J$.
Therefore, it suffices to prove that $Q$ is semisimple, or equivalently that $c_0 \ne 0$.
Suppose $c_0 = 0$, on the contrary.
The elements $z$ in $J \cap P(H, K)$ which is defined by the same formula as (\ref{eq:2.9})
is again a non-zero primitive with $z^p = 0$.
Here, to see that $z$ is a primitive, notice that the element $\gamma_z$ in $K \otimes K$ ($= k^{G \times G}$)
defined so as by (\ref{eq:2.5}) is zero since $\gamma_z^p = 0$.
We have thus proved $c_0 \ne 0$, as desired. \hfill $\qed$

\begin{remark} \label{rem:3.1}
  For the condition (iv), we have actually used the condition (c) in Proposition \ref{prop:1.3}.
  Therefore the proof above, applied to the special case $H = u(L)$, proves (c) $\iff$ (d) $\iff$ (e) in the proposition.
\end{remark}

\begin{remark} \label{rem:3.2}
  Suppose $\ch k = p > 0$, and that $k$ is algebraically closed.
  The coordinate Hopf algebra $\mathcal{O}(\boldsymbol{\alpha}_p)$ of the finite group scheme $\boldsymbol{\alpha}_p$
  is given by $k[z]/(z^p)$ with $z$ a primitive; this is selfdual.
  Demazure and Gabriel's result \cite[IV, \S 3, 3.7]{bib:DG} cited in the Introduction
  is translated into the language of Hopf algebras, as follows:

  {\it
    If $H$ is a commutative finitely generated Hopf algebra such that
    (a) $H$ does not contain any non-trivial idempotent, and
    (b) $H$ does not have any quotient Hopf algebra isomorphic to $\mathcal{O}(\boldsymbol{\alpha}_p)$,
    then $H \cong kG$, a group algebra of some (necessarily abelian) group $G$.
  }

  By using the same idea as Sweedler's \cite{bib:Sw2},
  we will deduce this result from our theorem, (iv) $\Longrightarrow$ (iii).
  It suffices to prove that every finite-dimensional subcoalgebra $C \subset H$ is spanned by grouplikes.
  For every $n \ge 0$, the ideal $I_n$ generated by $h^{p^n}$ ($h \in H^+$) is a Hopf ideal.
  One sees easily that the Hopf algebra
  $H / I_n$ is finite-dimensional and local.
  The dual $(H / I_n)^*$ is irreducible, and by the assumption (b),
  it does not contain any primitive $z \ne 0$ with $z^p = 0$.
  By our theorem, (iv) $\Longrightarrow$ (iii), $H / I_n$ is spanned by grouplikes.
  One sees $I_n \subset (H^+)^{p^n}$, and $\bigcap_{n} (H^+)^n = 0$
  by the assumption (a) and the Krull intersection theorem; see \cite[Theorem 2.10]{bib:Sw2}.
  Hence, $C \subset H / I_n$ for $n \gg 0$, which implies that $C$ is spanned by grouplikes, as desired.

  On the other hand,
  Sweedler's \cite[Theorem 4.1]{bib:Sw2} easily follows from that theorem specialized when $k = \bar{k}$,
  and the Hopf algebra in question is finitely generated; see also \cite[Sect. 5.7]{bib:Mo}.
  The specialized theorem is formulated just as above, except that the assumption (b) is replaced by
  \begin{itemize}
  \item [(b')] $H$ is cosemisimple.
  \end{itemize}
  For the thus formulated result, the proof above is valid since by (b'),
  $H / I_n$ is necessarily cosemisimple, as was proved by Sweedler \cite[Corollary 1.10]{bib:Sw2},
  whence it is spanned by grouplikes, by our theorem, (i) $\Longrightarrow$ (iii) (applied to $(H / I_n)^*$).
\end{remark}

\end{document}